\documentclass[11pt]{article}%
\usepackage[english]{babel}
\usepackage{indentfirst}
\usepackage{latexsym}
\usepackage[colorlinks=true, bookmarksnumbered=true, bookmarksopen=true,
bookmarksopenlevel=3, pdfstartview=FitH, linkcolor=cyan, pdfmenubar=true,
pdftoolbar=true, bookmarks=true,citecolor=cyan, urlcolor=magenta,
filecolor=cyan,menucolor=black,plainpages=false,pdfpagelabels]{hyperref}
\usepackage[lmargin=3.0cm,tmargin=2.5cm,bmargin=2.0cm,rmargin=2.7cm]%
{geometry}
\usepackage{amsbsy, amsmath,amsfonts,amssymb,amsthm}
\usepackage{times}
\usepackage{graphicx}
\usepackage{amsmath}
\usepackage{amsfonts}
\usepackage{amssymb}
\usepackage[english]{babel}
\usepackage[latin1]{inputenc}
\usepackage[T1]{fontenc}
\usepackage{amsmath,amssymb,graphics,amsthm}
\usepackage[usenames]{color}%
\numberwithin{equation}{section}
\newtheorem{prop}{Proposition}[section]

\newtheorem{teo}[prop]{Theorem}

\newtheorem{obs}[prop]{Remark}
\newtheorem{lema}[prop]{Lemma}

\newcommand\reallywidehat[1]{\savestack{\tmpbox}{\stretchto{  \scaleto{    \scalerel*[\widthof{\ensuremath{#1}}]{\kern-.6pt\bigwedge\kern-.6pt}    {\rule[-\textheight/2]{1ex}{\textheight}}  }{\textheight}}{0.5ex}}\stackon[1pt]{#1}{\tmpbox}}
\def\fin { \vskip 0pt \hfill $\diamond$ \vskip 12pt}
\begin{document}

\title{Global well-posedness for the fractional Boussinesq-Coriolis system with
stratif{}ication in a framework of Fourier-Besov type \\\vspace{0.3cm}}
\author{{Leithold L. Aurazo-Alvarez$^{1}$}, \ {Lucas C. F. Ferreira$^{2}$%
}{\thanks{{Corresponding author. }\newline{E-mail adresses:
aurazoall@gmail.com (L.L. Aurazo-Alvarez), lcff@ime.unicamp.br (L.C.F.
Ferreira).}\newline{L.L. Aurazo-Alvarez was financed by CAPES (Finance Code
001) and Cnpq (Grant 141166/2019-3), Brazil. \newline L.C.F. Ferreira was
supported by FAPESP and CNPq, Brazil.}}}\\{\vspace{-0.3cm}}\\{\small $^{1,2}$ University of Campinas, IMECC-Department of Mathematics} \\{\small CEP 13083-859, Campinas, SP, Brazil.}}
\date{}
\maketitle

\begin{abstract}
\medskip

We establish the global well-posedness of the 3D fractional
Boussinesq-Coriolis system with stratif{}ication in a framework of Fourier
type, namely spaces of Fourier-Besov type with underlying space being Morrey
spaces (FBM-spaces, for short). Under suitable conditions and rescaled density
f{}luctuation, the result is uniform with respect to the Coriolis and
stratification parameters. We cover the critical case of the dissipation,
namely half-Laplacian, in which the nonlocal dissipation has the same
differential order as the nonlinearity and balances critically the scaling of
the quadratic nonlinearities. As a byproduct, considering trivial initial
temperature and null stratification, we also obtain well-posed results in
FBM-spaces for the fractional Navier-Stokes-Coriolis system as well as for the
Navier-Stokes equations with critical dissipation. Moreover, since small
conditions are taken in the weak norm of FBM-spaces, we can consider some
initial data with arbitrarily large $H^{s}$-norms, $s\geq0.$\vspace{0.1cm}


 \vspace{0.1cm}

\end{abstract}

\renewcommand{\abstractname}{Abstract}

\section{Introduction}

We are concerned with the initial value problem (IVP) for the 3D fractional
Boussinesq-Coriolis equations with stratif{}ication (FBCS)
\begin{equation}
\,\,\,\left\{
\begin{split}
&  \partial_{t}u+\nu(-\Delta)^{\alpha}u+\Omega e_{3}\times u+(u\cdot
\nabla)u+\nabla p=g\theta e_{3},\\
&  \partial_{t}\theta+k(-\Delta)^{\alpha}\theta+(u\cdot\nabla)\theta
=-\mathcal{N}^{2}u_{3},\\
&  div\,u=0\,\,\mbox{for}\,\,(x,t)\in\mathbb{R}^{3}\times(0,\infty
)\,\,\mbox{and}\\
&  u(x,0)=u_{0}(x),\,\,\theta(0,x)=\theta_{0}(x)\,\,\mbox{for}\,\,x\in
\mathbb{R}^{3},
\end{split}
\right.  \label{FSBC}%
\end{equation}
where $u=(u_{1}(x,t),u_{2}(x,t),u_{3}(x,t))$, $\theta=\theta(x,t)$ and
$p=p(x,t)$ stand for the f{}luid velocity, the density f{}luctuation and the
pressure of the f{}luid, respectively. The kinetic viscosity, the thermal
diffusivity and the gravity are respectively represented by the positive
constants $\nu,k$ and $g$. The term $\Omega e_{3}\times u$ denotes the
so-called Coriolis force where the parameter $\Omega\neq0$ is the speed of
rotation of the fluid around the vertical unit vector $e_{3}=(0,0,1).$ The
term $g\theta e_{3}$ comes from the Boussinesq approximation (see
\cite{Cushman-Roisin}) in which density variations influence proportionally in
the gravitational term, while $\mathcal{N}^{2}u_{3}$ carries information about
the stratif{}ication effects where the stratif{}ication parameter
$\mathcal{N}>0$ represents the Brunt-V\"{a}is\"{a}l\"{a} wave frequency
related to the buoyancy of the f{}luid. Moreover, the divergence-free vector
$u_{0}(x)=(u_{0,1},u_{0,2},u_{0,3})$ is the initial velocity and the scalar
function $\theta_{0}=\theta_{0}(x)$ is the initial density disturbance.

The operator $(-\Delta)^{\alpha}$ is the fractional power of the minus
Laplacian and we consider the range $\frac{1}{2}\leq\alpha<\frac{5}{2}.$ In
the half-Laplacian case $\alpha=1/2$, we have the critical dissipation in the
sense that the nonlocal dissipation has the same differential order as the
nonlinearity and balances critically the scaling of the quadratic
nonlinearities in (\ref{FSBC}), which introduces further difficulties.

The two effects, one arising from the Coriolis force related to the parameter
$\Omega$ and the another from the stable stratif{}ication related to the
parameter $\mathcal{N}$, play an important role in large scale atmosphere
dynamics in which the density fluctuation $\theta$ is considered to depend
only on the potential temperature. They are also used in the study of
large-scale motions of the oceans. For further details, we refer the reader to
\cite{Cushman-Roisin, Chemin1, IMT2017}. In view of the analysis of asymptotic
regimes as $\Omega$ and $\mathcal{N}$ go to infinity (fast
oscillating/strongly stratified limit), an important subject is to find
frameworks in which system (\ref{FSBC}) is well-posed with bound of the
solution and initial-data size (or existence-time) independent of those
parameters, that is, uniform with respect to $\Omega$ and $\mathcal{N}$ (see,
e.g., \cite{BMN-4, BMN-2, GIMS2008}), at least uniformly for large $\Omega$
and $\mathcal{N}$ (say $\left\vert \Omega\right\vert ,\left\vert
\mathcal{N}\right\vert >c_{0}$ where $c_{0}$ is a positive universal
constant). Moreover, motivated by the study of statistical properties of
turbulence, spaces containing functions nondecaying at infinity are of
interest (see, e.g., \cite{Foias1, GIMS2008}). For these purposes, we employ a
framework of Fourier type introduced in \cite{Ferreira-Lima2014} to analyze
active scalar equations, namely Fourier-Besov-Morrey spaces $\mathcal{FN}%
_{q,\mu,r}^{s}$ (FBM-spaces, for short). These spaces are of Fourier-Besov
type with underlying space being Morrey spaces and, considering the same
scaling, are larger than classical Fourier-Besov spaces $\dot{FB}_{q,r}^{s}$
\cite{Iwabuchi-Takada2014} (see Section 2.1, p.4, for the definition and
details). Throughout this paper, spaces of scalar and vector functions are
denoted in the same way.

Considering $\theta\equiv0$, $\mathcal{N}=0$ and $\Omega=0$ in (\ref{FSBC}),
we have the 3D fractional Navier-Stokes equations (3DFNS) for which there is a
wide literature about existence of global mild solutions in different critical
frameworks. A Banach space $X$ is said to be critical for (3DFNS) if
$\left\Vert f(x)\right\Vert _{X}\approx\left\Vert \lambda^{2\alpha-1}f(\lambda
x)\right\Vert _{X}$ for all $\lambda>0,$ that is, the norm is invariant under
the scaling $f(x)\rightarrow\lambda^{2\alpha-1}f(\lambda x).$ Let us start by
briefly reviewing some results for the case $\alpha=1$ which corresponds to
the celebrated 3D Navier-Stokes equations. Without making a complete list, we
would like to mention the global well-posedness results with small initial
data in critical spaces such as Lebesgue $L^{n}$ \cite{Kato2}, Besov $\dot
{B}_{q,\infty}^{\frac{n}{q}-1}$ \cite{Cannone3}, Morrey $\mathcal{M}_{q,n-q}$
\cite{Kato-Morrey, G-M1989}, Fourier-Besov $\dot{FB}_{q,\infty}^{n-1-\frac
{n}{q}}$ \cite{Iwabuchi-Takada2014, Konieczny-Yoneda2011}, Besov-Morrey
$\dot{N}_{q,\mu,\infty}^{\frac{n-\mu}{q}-1}$ \cite{Kozo-Yamazaki1994,
Mazzucato}, and $BMO^{-1}$ \cite{Koch-Tataru2001}, among others. For further
details, we refer the reader to the book \cite{Lemarie2002}. Also, there are
some results about ill-posedness, for instance, see \cite{Bourgain1, Wang1}
and their references for results in the space $\dot{B}_{\infty,r}%
^{-1}(\mathbb{R}^{3})$ with $1\leq r\leq\infty$.

Still for (3DFNS), but with the fractional dissipation $\alpha$, the author of
\cite{Lions1969} obtained the global existence of classical solutions for
$\alpha\geq\frac{5}{4}$, without any smallness condition. Nevertheless, the
global well-posedness in the case $\alpha<\frac{5}{4}$ is more subtle and an
outstanding open problem. For small initial data in critical spaces, there are
global well-posedness results in the Besov space $\dot{B}_{q,r}^{1-2\alpha
+\frac{3}{q}}(\mathbb{R}^{3})$ \cite{Wu2005}, in the largest critical space
$\dot{B}_{\infty,\infty}^{1-2\alpha}(\mathbb{R}^{3})$ with $\frac{1}{2}%
<\alpha<1$ \cite{Yu-Zhai2012}, in the Triebel-Lizorkin space $\dot{F}%
_{\frac{3}{\alpha-1},2}^{-\alpha}(\mathbb{R}^{3})$ with $1<\alpha<\frac{5}{4}$
\cite{Deng-Yao2014}, in the Fourier-Besov space $\dot{FB}_{q,r}^{4-2\alpha
-\frac{3}{q}}(\mathbb{R}^{3})$ for $1\leq q\leq r\leq2$ and $\frac{1}%
{2}<\alpha\leq\frac{5}{2}-\frac{3}{2q}$ and for $\alpha=\frac{1}{2}$ with
$r=1$ and $1\leq q\leq\infty$ \cite{Xiao-CW2014}, and in the
Fourier-Besov-Morrey space $\mathcal{FN}_{q,\mu,r}^{4-2\alpha-\frac{3-\mu}{q}%
}(\mathbb{R}^{3})$ for $\frac{1}{2}<\alpha<\frac{5}{2}-\frac{3-\mu}{2q}$,
$0\leq\mu<3$, $1\leq q<\infty$ and $1\leq r\leq\infty$ (see
\cite{Bara-Toumli2019} with $\Omega=0$). Moreover, we mention ill-posedness
results in the largest critical space $\dot{B}_{\infty,\infty}^{1-2\alpha
}(\mathbb{R}^{3})$ for $1\leq\alpha<\frac{5}{4}$ \cite{CS2012} and in the
Triebel-Lizorkin space $\dot{F}_{\frac{3}{\alpha-1},r}^{-\alpha}%
(\mathbb{R}^{3})$ for $r>2$ and $1<\alpha<\frac{5}{4}$ \cite{Deng-Yao2014}.

Another relevant model covered by (\ref{FSBC}) corresponds to the case
$\theta\equiv0$, $\mathcal{N}=0$ and general $\Omega$, namely the fractional
Navier-Stokes-Coriolis system. For the value $\alpha=1$, that is, the
classical Navier-Stokes-Coriolis system, we have $\Omega$-uniform global
well-posedness results for small initial data in critical Fourier
transform-based functional spaces, for clarity, critical with respect to the
(3DFNS)-scaling. For instance, Hieber and Shibata \cite{Hieber2010} showed
existence of a unique global mild solution in $H^{\frac{1}{2}}(\mathbb{R}%
^{3})$, where the smallness condition is uniform w.r.t. $\Omega$. After, Giga
et al \cite{GIMS2008} proved $\Omega$-uniform global well-posedness in
$FM_{0}^{-1}(\mathbb{R}^{3})$ which can be identified with $\dot{FB}%
_{1,1}^{-1}$ and permits to consider spatially nondecaying and almost periodic
initial-data. The $\Omega$-uniform well-posedness in $\dot{FB}_{q,\infty
}^{2-\frac{3}{q}}(\mathbb{R}^{3})$ with $1<q\leq\infty$ and $\dot{FB}%
_{1,2}^{-1}(\mathbb{R}^{3})$ were proved by Konieczny and Yoneda
\cite{Konieczny-Yoneda2011} and Iwabuchi and Takada \cite{Iwabuchi-Takada2014}%
, respectively. Employing the framework of FBM-spaces, Almeida \textit{et al.}
\cite{AFL2017} obtained the $\Omega$-uniform global well-posedness in
$\mathcal{FN}_{q,\mu,\infty}^{2-\frac{3-\mu}{q}}(\mathbb{R}^{3})$ where $1\leq
q<\infty$ and $0\leq\mu<3$ with $\mu\neq0$ when $q=1$. Moreover, we have
ill-posedness in $\dot{FB}_{1,r}^{-1}(\mathbb{R}^{3})$ when $2<r\leq\infty$
(see \cite{Iwabuchi-Takada2014}). For general index $\alpha$, we have $\Omega
$-uniform global well-posedness results in the Lei-Lin-type space
$\mathcal{X}^{1-2\alpha}(\mathbb{R}^{3})$ with $\frac{1}{2}\leq\alpha\leq1$
\cite{Wang-Wu1}, in the Fourier-Besov space $\dot{FB}_{q,r}^{4-2\alpha
-\frac{3}{q}}(\mathbb{R}^{3})$ with $\frac{2}{3}<\alpha<\frac{5}{3}-\frac
{1}{q}$, $2\leq q\leq\infty$ and $1\leq r\leq\infty$ \cite{Wang-Wu2}, and in
the FBM-space $\mathcal{FN}_{q,\mu,r}^{4-2\alpha-\frac{3-\mu}{q}}%
(\mathbb{R}^{3})$ with $\frac{1}{2}<\alpha\leq\frac{5}{2}-\frac{3-\mu}{2q}$,
$0\leq\mu<3$, $1\leq q<\infty$ and $1\leq r\leq2$ \cite{Bara-Toumli2019}.
\newline

For (\ref{FSBC}) with $\alpha=1,$ Sun and Cui \cite{Sun-Cui2019} showed the
global well-posedness with small initial-data in the critical space $\dot
{FB}_{q,r}^{2-\frac{3}{q}}(\mathbb{R}^{3})$ for $1<q\leq\infty$ and $1\leq
r<\infty$ and in $\dot{FB}_{1,r}^{-1}(\mathbb{R}^{3})$ for $1\leq r\leq2$.
Moreover, they also proved the ill-posedness in the Fourier-Besov space
$\dot{FB}_{1,r}^{-1}(\mathbb{R}^{3})$ for $2<r\leq\infty$.

Finally, it is worth to mention that another class of study investigates the
dispersive effect of the associated linear semigroup in order to show global
well-posedness, without smallness conditions, for geophysical (and related)
models with rotation and/or stratif{}ication sufficiently large in comparison
with the initial-data norm (see, e.g., \cite{BMN-4, Chemin1,
Takada2013,IMT2017, BMN-3, Charve2, KMY2012}). More precisely, they show
global well-posedness and asymptotic results in $H^{s}$ for large-enough
values of the parameters $\Omega$ and/or $\mathcal{N}$, where the bounds from
below for them depend on the initial-data size. See also
\cite{Angulo-Ferreira} for an extension of these ideas to obtain global
well-posedness in the context of Besov spaces for the case $\theta\equiv0$,
$\mathcal{N}=0$ and $\alpha=1$ in (\ref{FSBC}) (Navier-Stokes-Coriolis equations).

In view of previous references, even in the classical dissipation case
$\alpha=1$, there are much less well-posedness results for system (\ref{FSBC})
than Navier-Stokes and Navier-Stokes-Coriolis equations. Motivated by that and
bearing in mind the papers \cite{AFL2017, Ferreira-Lima2014}, we consider the
framework of FBM-spaces and obtain a global well-posedness result for
(\ref{FSBC}) (see Theorem \ref{TH1}), which provides a larger critical
initial-data class for global well-posedness that allows to consider singular
and nondecaying initial data. Since small conditions are taken in the weak
norm of FBM-spaces, we can consider some initial data with arbitrarily large
$H^{s}$-norms, $s\geq0,$ e.g., suitable cut-offs in Fourier variables of
homogeneous functions (see Lemma \ref{Lem-Properties} (v)). Choosing suitable
indexes of the spaces, we cover the fractional dissipation range $\frac{1}%
{2}\leq\alpha<\frac{5}{2}$ including the critical dissipation $\alpha=1/2.$
For that, besides the norm of the persistence space, we employ an auxiliary
norm of Chemin-Lerner type based on $\mathcal{FN}_{q,\mu,r}^{s}$-spaces in
order to estimate the linear and bilinear terms in (\ref{FSBC}). Considering
the rescaled variable $v=(u,\sqrt{g}\theta/\mathcal{N})$ and assuming, for
instance, $\mathcal{N}\sqrt{g}/2\leq\left\vert \Omega\right\vert
\leq2\mathcal{N}\sqrt{g},$ the well-posedness result is uniform with respect
to the parameters $\mathcal{N}$ and $\Omega$. The present paper is part of
Ph.D thesis \cite{Aurazo}. After concluding this work, we became aware of the
preprint \cite{AAO} in which the authors independently proved a well-posedness
result with $\alpha=1$ similar to the one obtained here.

The outline of this manuscript is as follows. In Section \ref{Preliminar}, we
first recall some definitions, notations, basic tools and the def{}inition of
Fourier-Besov-Morrey spaces. After, we define the fractional
Boussinesq-Coriolis-Stratif{}ication semigroup, suitable time-dependent
spaces, and mild solutions for (\ref{FSBC}). In Section \ref{main} we state
our main result and make some comments about it. Section \ref{Proof} is
devoted to prove linear and bilinear estimates in our functional setting as
well as the proof of the main result.

\section{Preliminaries}

\label{Preliminar} In this section we recall some basic properties about
Fourier-Besov-Morrey spaces, the fractional Boussinesq-Coriolis-Stratif{}%
ication semigroup and other analysis tools we shall use throughout this work.

\subsection{Fourier-Besov-Morrey spaces}

Fourier-Besov-Morrey spaces (FBM-spaces) are constructed by mean of a kind of
localization procedure in Fourier variables on the well-known Morrey spaces
(see \cite{Ferreira-Lima2014}). Morrey spaces are def{}ined as follows. Let
$B_{d}(x_{0})$ be the open ball in $\mathbb{R}^{n}$ centered at $x_{0}$ and
with radius $d>0$. For $1\leq q<\infty$ and $0\leq\mu<n$, the homogeneous
Morrey space $\mathcal{M}_{q,\mu}=\mathcal{M}_{q,\mu}(\mathbb{R}^{n})$ is the
space of all $f\in L_{loc}^{q}$ such that
\begin{equation}
\parallel f\parallel_{q,\mu}=\displaystyle{\sup_{x_{0}\in\mathbb{R}^{n},d>0}%
}d^{-\frac{\mu}{q}}\parallel f\parallel_{L^{q}(B_{d}(x_{0}))}<\infty.
\label{Mspace}%
\end{equation}
In the case $q=1$, $\mathcal{M}_{1,\mu}$ is a subspace of Radon measures and
the $L^{1}$-norm in (\ref{Mspace}) should be understood as the total variation
of the measure $f$ on $B_{d}(x_{0})$. The space $\mathcal{M}_{q,\mu}$ endowed
with $\parallel\cdot\parallel_{q,\mu}$is a Banach space. For more details, we
refer the reader to \cite{Kato-Morrey} and their references. \newline\qquad
Next we recall some notations, definitions and properties about
Littlewood-Paley decomposition. First recall the notations $\mathcal{S}%
(\mathbb{R}^{n})$ and $\mathcal{S}^{\prime}(\mathbb{R}^{n})$ to stand for the
Schwartz class and the space of tempered distributions, respectively. Also,
the Fourier transform of $u$ in $\mathcal{S}^{\prime}(\mathbb{R}^{n})$ is
denoted by $\hat{u}$.

Consider the ring $\mathcal{C}=\{\xi\in\mathbb{R}^{n};\frac{3}{4}\leq\mid
\xi\mid\leq\frac{8}{3}\}$ and $\varphi$ a smooth function supported in
$\mathcal{C}$ satisfying $0\leq\varphi\leq1$ and
\begin{equation}
\displaystyle{\sum_{j\in\mathbb{Z}}}\varphi_{j}(\xi
)=1,\,\,\mbox{for all}\,\,\xi\neq0,\,\mbox{where}\,\,\varphi_{j}(\xi
)=\varphi(2^{-j}\xi). \label{somatoria}%
\end{equation}
For each $j\in\mathbb{Z}$, the localization operator $\Delta_{j}$ is defined
via Fourier transform as
\[
\lbrack\Delta_{j}f]\symbol{94}(\xi)=\varphi_{j}(\xi)\hat{f}(\xi)\text{
in}\,\,\,\mathcal{S}^{\prime}(\mathbb{R}^{n}),
\]
and the operator $S_{j}$ as
\[
S_{j}f=\displaystyle{\sum_{k\leq j-1}}\Delta_{k}f\,\,\text{\ in}%
\,\,\,\mathcal{S}^{\prime}(\mathbb{R}^{n}).
\]
Since $supp\,\varphi_{j}\subset2^{j}\mathcal{C}$ we verify that
\[
\Delta_{j}\Delta_{k}f=0,\,\,\mbox{if}\,\,\mid j-k\mid\geq2,\,\,
\]
and
\[
\mbox{and}\,\,\Delta_{j}[S_{k-1}f\Delta_{k}g]=0,\,\,\mbox{if}\,\,\mid
j-k\mid\geq5.
\]
Let $1\leq q<\infty$, $0\leq\mu<n$, $1\leq r<\infty$ and $s\in\mathbb{R}$. The
Fourier-Besov-Morrey space $\mathcal{FN}_{q,\mu,r}^{s}$ is the set of all
distributions $f\in\mathcal{S}^{\prime}/\mathcal{P}$, where $\mathcal{P}$ is
the set of all polynomials in $\mathbb{R}^{n}$, such that $\varphi_{j}\hat
{f}\in\mathcal{M}_{q,\mu}$, for all $j\in\mathbb{Z}$, and
\begin{equation}
\parallel f\parallel_{\mathcal{FN}_{q,\mu,r}^{s}}=\left\{
\begin{split}
&  \displaystyle{(\displaystyle{\sum_{j\in\mathbb{Z}}}\displaystyle{(2^{js}%
\parallel\varphi_{j}\hat{f}\parallel_{q,\mu})^{r}})}^{1/r}<\infty,\,\,1\leq
r<\infty\\
&  \displaystyle{\sup_{j\in\mathbb{Z}}}\,2^{js}\parallel\varphi_{j}\hat
{f}\parallel_{q,\mu}<\infty,\,\,r=\infty.
\end{split}
\right.  \label{FBM}%
\end{equation}
The pair $(\mathcal{FN}_{q,\mu,r}^{s},\parallel\cdot\parallel_{\mathcal{FN}%
_{q,\mu,r}^{s}})$ is a Banach space. The next lemma contains some basic
properties of Morrey (see \cite{Kato-Morrey}) and FBM-spaces (see
\cite{Ferreira-Lima2014} and \cite{AFL2017}).

\begin{lema}
\label{Lem-Properties}Let $s_{1},s_{2}\in\mathbb{R}$, $1\leq p_{1},p_{2}%
,p_{3}<\infty$ and $0\leq\mu_{1},\mu_{2},\mu_{3}<n$.

\begin{itemize}
\item[(i)] (H\"older's inequality) Let $\frac{1}{p_{3}}=\frac{1}{p_{2}}%
+\frac{1}{p_{1}}$ and $\frac{\mu_{3}}{p_{3}}=\frac{\mu_{2}}{p_{2}}+\frac
{\mu_{1}}{p_{1}}$. If $f_{i} \in\mathcal{M}_{p_{i},\mu_{i}}$ for $i=1,2$, then
$f_{1}f_{2} \in\mathcal{M}_{p_{3},\mu_{3}}$ and
\[
\parallel f_{1}f_{2}\parallel_{p_{3},\mu_{3}} \leq\parallel f_{1}%
\parallel_{p_{1},\mu_{1}} \parallel f_{2}\parallel_{p_{2},\mu_{2}}.
\]

\item[(ii)] (Young's inequality) If $\varphi\in L^{1}(\mathbb{R}^{n})$ and
$g\in\mathcal{M}_{ p_{1},\mu_{1}}$ then
\[
\parallel\varphi\ast g\parallel_{p_{1},\mu_{1}} \leq\parallel\varphi
\parallel_{1} \parallel g\parallel_{p_{1},\mu_{1}},
\]
where $\ast$ denotes the standard convolution operator.

\item[(iii)] (Bernstein-type inequality) Let $p_{2}\leq p_{1}$ be such that
$\frac{n-\mu_{1}}{p_{1}}\leq\frac{n-\mu_{2}}{p_{2}}$. If $A>0$ and
$supp\,(\hat{f})\subset\{\xi\in\mathbb{R}^{n};\mid\xi\mid\leq A2^{j}\}$, then
\begin{equation}
\parallel\xi^{\beta}\hat{f}\parallel_{p_{2},\mu_{2}}\leq C2^{j\mid\beta
\mid+j(\frac{n-\mu_{2}}{p_{2}}-\frac{n-\mu_{1}}{p_{1}})}\parallel\hat
{f}\parallel_{p_{1},\mu_{1}},
\end{equation}
where $\beta$ is the multi-index, $j\in\mathbb{Z}$, and $C>0$ is a constant
independent of $j,\xi$ and $f$.

\item[(iv)] (Sobolev-type embedding) For $p_{2}\leq p_{1}$ and $s_{2}\leq
s_{1}$ satisfying $s_{2}+\frac{n-\mu_{2}}{p_{2}}=s_{1}+\frac{n-\mu_{1}}{p_{1}%
}$, we have the continuous inclusion
\[
\mathcal{FN}_{p_{1},\mu_{1},r_{1}}^{s_{1}}\subset\mathcal{FN}_{p_{2},\mu
_{2},r_{2}}^{s_{2}},
\]
for all $1\leq r_{1}\leq r_{2}\leq\infty$.

\item[(v)] The space $\mathcal{FN}_{q,\mu,r}^{s}$ contains homogeneous
functions of degree $-k=s-n+\frac{n-\mu}{q}$.
\end{itemize}
\end{lema}

\subsection{Fractional Boussinesq-Coriolis-Stratif{}ication semigroup and mild
solutions}

In this subsection, by following \cite{KMY2012} (see also \cite{Sun-Cui2019}
and \cite{IMT2017}), we recall how to rewrite system (\ref{FSBC}) as an
integral equation. Considering $N=\mathcal{N}\sqrt{g}$, $v=(v^{1},v^{2}%
,v^{3},v^{4})=(u^{1},u^{2},u^{3},\sqrt{g}\theta/\mathcal{N})$, $v_{0}%
=(v_{0}^{1},v_{0}^{2},v_{0}^{3},v_{0}^{4})=(u_{0}^{1},u_{0}^{2},u_{0}%
^{3},\sqrt{g}\theta_{0}/\mathcal{N})$, and $\tilde{\nabla}=(\partial
_{1},\partial_{2},\partial_{3},0)$, we can convert system (\ref{FSBC}) to
\begin{equation}
\left\{
\begin{split}
&  \partial_{t}v+\mathcal{A}v+\mathcal{B}v+\tilde{\nabla}p=-(v\cdot
\tilde{\nabla})v,\,\,\mbox{in}\,\,\mathbb{R}^{3}\times(0,\infty)\\
&  \tilde{\nabla}\cdot v=0,\,\,\mbox{in}\,\,\mathbb{R}^{3}\times(0,\infty)\\
&  v(x,0)=v_{0}(x),\,\,\mbox{in}\,\,\mathbb{R}^{3},
\end{split}
\right.  \label{SFBQ2}%
\end{equation}
where
\[
\mathcal{A}=\left(
\begin{matrix}
\nu(-\Delta)^{\alpha} & 0 & 0 & 0\\
0 & \nu(-\Delta)^{\alpha} & 0 & 0\\
0 & 0 & \nu(-\Delta)^{\alpha} & 0\\
0 & 0 & 0 & k(-\Delta)^{\alpha}%
\end{matrix}
\right)  \,\,\,\,\mbox{and}\,\,\,\,\mathcal{B}=\left(
\begin{matrix}
0 & -\Omega & 0 & 0\\
\Omega & 0 & 0 & 0\\
0 & 0 & 0 & -N\\
0 & 0 & N & 0
\end{matrix}
\right)  .
\]
The fractional Boussinesq-Coriolis-Stratif{}ication semigroup
((FBCS)-semigroup), namely the semigroup associated to the linear part of
(\ref{SFBQ2}), is denoted by $\{S_{\Omega,N}^{\alpha}(t)\}_{t\geq0}.$ For
$\nu=k$, we have the following expression for $\{S_{\Omega,N}^{\alpha
}(t)\}_{t\geq0}$ in Fourier variables
\begin{equation}
\widehat{S_{\Omega,N}^{\alpha}(t)f}=\mathcal{F}^{-1}[\cos(\frac{\mid\xi
\mid^{\prime}}{\mid\xi\mid}t)e^{-\nu t\mid\xi\mid^{2\alpha}}M_{1}(\xi)\hat
{f}+\sin(\frac{\mid\xi\mid^{\prime}}{\mid\xi\mid}t)e^{-\nu t\mid\xi
\mid^{2\alpha}}M_{2}(\xi)\hat{f}+e^{-\nu t\mid\xi\mid^{2\alpha}}M_{3}(\xi
)\hat{f}], \label{Semigroup}%
\end{equation}
where, for each $\xi=(\xi_{1},\xi_{2},\xi_{3})\in\mathbb{R}^{3},$
\[
\mid\xi\mid=\sqrt{\xi_{1}^{2}+\xi_{2}^{2}+\xi_{3}^{2}},\,\,\,\mid\xi
\mid^{\prime}=\sqrt{N^{2}\xi_{1}^{2}+N^{2}\xi_{2}^{2}+\Omega^{2}\xi_{3}^{2}},
\]%
\[
M_{1}(\xi)=\frac{1}{\mid\xi\mid^{\prime}{}^{2}}\left(
\begin{matrix}
\Omega^{2}\xi_{3}^{2} & 0 & -N^{2}\xi_{1}\xi_{3} & \Omega N\xi_{2}\xi_{3}\\
0 & \Omega^{2}\xi_{3}^{2} & -N^{2}\xi_{2}\xi_{3} & -\Omega N\xi_{1}\xi_{3}\\
-\Omega^{2}\xi_{1}\xi_{3} & -\Omega^{2}\xi_{2}\xi_{3} & N^{2}(\xi_{1}^{2}%
+\xi_{2}^{2}) & 0\\
\Omega N\xi_{2}\xi_{3} & -\Omega N\xi_{1}\xi_{3} & 0 & N^{2}(\xi_{1}^{2}%
+\xi_{2}^{2})
\end{matrix}
\right)  ,
\]

\[
M_{2}(\xi)=\frac{1}{\mid\xi\mid\mid\xi\mid^{\prime}}\left(
\begin{matrix}
0 & -\Omega\xi_{3}^{2} & \Omega\xi_{2}\xi_{3} & N\xi_{1}\xi_{2}\\
\Omega\xi_{3}^{2} & 0 & -\Omega\xi_{1}\xi_{3} & N\xi_{2}\xi_{3}\\
-\Omega\xi_{2}\xi_{3} & \Omega\xi_{1}\xi_{3} & 0 & -N(\xi_{1}^{2}+\xi_{2}%
^{2})\\
-N\xi_{1}\xi_{3} & -N\xi_{2}\xi_{3} & N(\xi_{1}^{2}+\xi_{2}^{2}) & 0
\end{matrix}
\right)  ,
\]
and
\[
M_{3}(\xi)=\frac{1}{\mid\xi\mid^{\prime}{}^{2}}\left(
\begin{matrix}
N^{2}\xi_{2}^{2} & -N^{2}\xi_{1}\xi_{2} & 0 & -\Omega N\xi_{2}\xi_{3}\\
-N^{2}\xi_{1}\xi_{2} & N\xi_{1}^{2} & 0 & \Omega N\xi_{1}\xi_{3}\\
0 & 0 & 0 & 0\\
-\Omega N\xi_{2}\xi_{3} & \Omega N\xi_{1}\xi_{3} & 0 & \Omega^{2}\xi_{3}^{2}%
\end{matrix}
\right)  .
\]
Thus, denoting the symbol of $S_{\Omega,N}^{\alpha}(t)$ in (\ref{Semigroup})
by $[S_{\Omega,N}^{\alpha}(t)]\symbol{94}(\xi),$ we can write
\[
\lbrack S_{\Omega,N}^{\alpha}(t)f]\symbol{94}(\xi)=[S_{\Omega,N}^{\alpha
}(t)]\symbol{94}(\xi)\hat{f}(\xi).
\]
\newline As pointed out by Sun and Cui \cite{Sun-Cui2019}, the components
$M_{j,k}^{l}$ of the matrix $M_{l}(\xi)$ can be estimated as
\begin{equation}
\mid M_{j,k}^{l}(\xi)\mid\leq L=:\max\{2,\frac{\mid\Omega\mid}{\mathcal{N}%
\sqrt{g}},\frac{\mathcal{N}\sqrt{g}}{\mid\Omega\mid}\}\,\,\mbox{for}\,\,1\leq
j,k\leq4,1\leq l\leq3\,\,\mbox{and for all}\,\,\xi\in\mathbb{R}^{3}.
\label{EstMatrixM}%
\end{equation}
\newline This bound will be used often in the estimates for the
(FBCS)-semigroup $\{S_{\Omega,N}^{\alpha}(t)\}_{t>0}$. \newline

Now we recall the extended Helmholtz projection operator $\tilde{\mathbb{P}%
}=(\tilde{\mathbb{P}}_{jk})_{4\times4}$
\begin{equation}
\tilde{\mathbb{P}}_{jk}=\left\{
\begin{split}
&  \delta_{jk}+R_{j}R_{k},\,\,1\leq j,k\leq3\\
&  \delta_{jk},\,\,\,\,\mbox{otherwise},
\end{split}
\right.  , \label{helm1}%
\end{equation}
where $\delta_{jk}$ denotes the Kronecker delta and $R_{j}$ is the Riesz
transform on$\mathbb{\ R}^{3}$ for each $j=1,2,3$. After applying this
operator on (\ref{SFBQ2}), we can employ Duhamel's principle in order to
formally convert system (\ref{SFBQ2}) to the integral equation
\begin{equation}
v(t)=S_{\Omega,N}^{\alpha}(t)v_{0}+B(v,v), \label{mildform}%
\end{equation}
where the bilinear operator $B(u,v)$ is def{}ined via the Fourier transform
as
\begin{equation}
\lbrack B(v,w)]\symbol{94}(\xi,t)=-\displaystyle{\int_{0}^{t}}[S_{\Omega
,N}^{\alpha}(t-\tau)]\symbol{94}(\xi)[\tilde{\mathbb{P}}]\symbol{94}(\xi
)(i\xi,0)^{T}\cdot\lbrack v\otimes w]\symbol{94}(\xi,\tau)d\,\tau.
\label{bilinearterm}%
\end{equation}
Observe that the projection $\tilde{\mathbb{P}}$ has the symbol
\[
\lbrack\tilde{\mathbb{P}}_{mk}]\symbol{94}(\xi)=\left\{
\begin{split}
&  \delta_{mk}-(\xi_{m}\xi_{k})/\left\vert \xi\right\vert ^{2},\,\,1\leq
m,k\leq3\\
&  \delta_{mk},\,\,\mbox{otherwise}
\end{split}
\right.  ,
\]
and the elements of the matrix $[v\otimes w]\symbol{94}(\xi,\cdot)$ are given
by
\[
([v\otimes w]\symbol{94})_{mk}=\hat{v}_{m}\ast\hat{w}_{k}(\xi,\cdot
),\,\,\mbox{for}\,\,1\leq m,k\leq4.
\]
Throughout the paper, solutions for the integral equation (\ref{mildform}) are
called mild solutions for system (\ref{FSBC}).

\subsection{Paraproduct formula, abstract f{}ixed point lemma and
time-dependent spaces}

In this subsection we recall Bony's paraproduct formula and an abstract
f{}ixed point lemma, and we def{}ine suitable time-dependent spaces based on
$\mathcal{FN}_{q,\mu,r}^{s}$ in which we shall look for solutions. \newline

For $f,g\in\mathcal{S}^{\prime}(\mathbb{R}^{n}),$ we can define Bony's
paraproduct operator $T_{f}$ \ as\
\begin{equation}
T_{f}(g)=\displaystyle{\sum_{j\in\mathbb{Z}}}S_{j-1}f\Delta_{j}g, \label{Top}%
\end{equation}
and the operator $R$
\begin{equation}
R(f,g)=\displaystyle{\sum_{j\in\mathbb{Z}}}\Delta_{j}f\tilde{\Delta}%
_{j}g\,\,\mbox{with}\,\,\tilde{\Delta}_{j}g=\displaystyle{\sum_{\mid
j^{\prime}-j\mid\leq1}}\Delta_{j^{\prime}}g. \label{Top2}%
\end{equation}
The symbol of the operator $\tilde{\Delta_{j}}$ in (\ref{Top2}) is denoted by
$\tilde{\varphi}_{j}$. Using these operators, we can write the commutative
product of $f$ by $g$ as
\begin{equation}
fg=T_{f}(g)+T_{g}(f)+R(f,g). \label{paraproduct}%
\end{equation}
Expression (\ref{paraproduct}) is known as Bony's paraproduct formula (see
\cite{Lemarie2002}). \newline

In order to avoid extensive fixed point computations, we shall use the
following lemma (see \cite{Lemarie2002}) for an abstract quadratic equation,
after obtaining the needed estimates.

\begin{lema}
\label{fixedpoint} Let $(X,\parallel\cdot\parallel)$ be a Banach space and
$B:X\times X\rightarrow$ X a bilinear operator satisfying $\parallel
B(x_{1},x_{2})\parallel\leq K\parallel x_{1}\parallel\parallel x_{2}\parallel$
for all $x_{1},x_{2}$, where $K>0$ is a constant. If $0<\varepsilon<\frac
{1}{4K}$ and $\parallel y\parallel\leq\varepsilon$, then the equation
$x=y+B(x,x)$ has a solution in $X$. Moreover, this solution is unique in the
closed ball $\{x\in X;\parallel x\parallel\leq2\varepsilon\}$ and $\parallel
x\parallel\leq2\parallel y\parallel$. The solution depends continuously on
$y$; in fact, if $\parallel\tilde{y}\parallel\leq\varepsilon$, $\tilde
{x}=\tilde{y}+B(\tilde{x},\tilde{x})$ and $\parallel\tilde{x}\parallel
\leq2\varepsilon$, then
\begin{equation}
\parallel x-\tilde{x}\parallel\leq(1-4K\varepsilon)^{-1}\parallel y-\tilde
{y}\parallel. \label{fixedpointdependence}%
\end{equation}

\end{lema}

Observe that, by writing $y=S_{\Omega,N}^{\alpha}(t)v_{0}$ and $B(v,w)$ as in
(\ref{bilinearterm}), the integral equation (\ref{mildform}) presents the form
$v=y+B(v,v)$ required in the above lemma.

Finally we def{}ine two time-dependent spaces. Let $1\leq p\leq\infty$,
$0<T\leq\infty$ and $I=(0,T)$. The Banach spaces $L^{p}(I;\mathcal{FN}%
_{q,\mu,r}^{s})$ and $\mathcal{L}^{p}(I;\mathcal{FN}_{q,\mu,r}^{s})$ are the
set of Bochner measurable functions from $I$ to $\mathcal{FN}_{q,\mu,r}^{s}$
with respective norms given by
\begin{equation}
\parallel f\parallel_{L^{p}(I;\mathcal{FN}_{q,\mu,r}^{s})}%
=\displaystyle{\parallel\parallel f(\cdot,t)\parallel_{\mathcal{FN}_{q,\mu
,r}^{s}}\parallel_{L^{p}(I)}}=\displaystyle{\parallel
\displaystyle{\ (\displaystyle{\sum_{j\in\mathbb{Z}}}(2^{js}\parallel
\varphi_{j}\hat{f}\parallel_{q,\mu})^{r})^{1/r}\parallel_{L^{p}(I)}}}%
\end{equation}
and
\begin{equation}
\parallel f\parallel_{\mathcal{L}^{p}(I;\mathcal{FN}_{q,\mu,r}^{s}%
)}=\displaystyle{\ (\displaystyle{\ \sum_{j\in\mathbb{Z}}(2^{js}%
\parallel\varphi_{j}\hat{f}\parallel_{L^{p}(I;\mathcal{M}_{q,\mu})})^{r}%
})^{1/r}}=\displaystyle{\ (\displaystyle{\ \sum_{j\in\mathbb{Z}}%
(2^{js}\parallel\parallel\varphi_{j}\hat{f}\parallel_{q,\mu}\parallel
_{L^{p}(I)})})^{r})^{1/r}}.
\end{equation}
We use the notation
\[
\mathcal{X}_{r}=:\mathcal{L}^{\infty}(I;\mathcal{FN}_{q,\mu,r}^{s}%
)\cap\mathcal{L}^{1}(I;\mathcal{FN}_{q,\mu,r}^{s+2\alpha}).
\]

\section{Main Result}

\label{main} In this section we state the main result of this work.

\begin{teo}
\label{TH1} Let $I=(0,\infty)$, $0\leq\mu< 3$, $1\leq q<\infty$ and
$s=4-2\alpha-\frac{3-\mu}{q}$.

Let $(\Omega,\mathcal{N})\in(\mathbb{R}-\{0\})^{2}$ and $v_{0}=(u_{0},\sqrt
{g}\theta_{0}/\mathcal{N})\in\mathcal{FN}_{q,\mu,r}^{s}$ be such that $u_{0}$
is a divergence-free vector f{}ield and $\theta_{0}$ a scalar f{}ield. Assume
that $\alpha$ and $r$ satisfy either the following conditions

\begin{itemize}
\item[(i)] $\frac{1}{2}<\alpha<\frac{5}{2}-\frac{3-\mu}{2q}$ and $1\leq
r\leq\infty$, or

\item[(ii)] $\frac{1}{2}<\alpha= \frac{5}{2}-\frac{3-\mu}{2q}$, $\mu=0$ and
$1\leq q\leq r\leq2$, or

\item[(iii)] $\alpha=\frac{1}{2}$ and $r=1$.
\end{itemize}

Then, there are two constants $\varepsilon=\varepsilon(\nu,\Omega
,\mathcal{N})>0$ and $C=C(\nu,\Omega,\mathcal{N})>0$ such that if $\parallel
v_{0}\parallel_{\mathcal{FN}_{q,\mu,r}^{s}}\leq\varepsilon$ then (\ref{FSBC})
has a unique global mild solution $v=(u,\sqrt{g}\theta/\mathcal{N}%
)\in\mathcal{X}_{r}$ such that $u$ is divergence-free and $v$ is unique in the
closed ball $\{v\in\mathcal{X}_{r};\parallel v\parallel_{X_{r}}\leq
C\varepsilon\}$. Moreover, the solution $(u,\sqrt{g}\theta/\mathcal{N})$ is
weakly time-continuous from $[0,\infty)$ to $\mathcal{S}^{\prime}%
(\mathbb{R}^{3})$ and depends continuously on the initial data $(u_{0}%
,\sqrt{g}\theta_{0}/\mathcal{N})$.

Assuming further that $\mathcal{N}\sqrt{g}/2\leq\mid\Omega\mid\leq
2\mathcal{N}\sqrt{g}$, we can take the constants $\varepsilon$ and $C$
independent of $\Omega$ and $\mathcal{N}$. So, for the rescaled unknown
variable $v=(u,\sqrt{g}\theta/\mathcal{N}),$ we obtain the uniform global
well-posedness with respect to both $\Omega$ and $\mathcal{N}$.
\end{teo}

\begin{obs}
Similar computations, using the fractional Stokes-Coriolis semigroup
$S_{\Omega}^{\alpha}(t)$, instead of the fractional
Boussinesq-Coriolis-Stratif{}ication semigroup, which is def{}ined via the
Fourier transform as
\[
\lbrack S_{\Omega}^{\alpha}(t)f]\symbol{94}(\xi)=e^{-\nu t\mid\xi\mid
^{2\alpha}}\cos(\frac{\Omega\xi_{3}}{\mid\xi\mid}t)\mathcal{I}+e^{-\nu
t\mid\xi\mid^{2\alpha}}\sin(\frac{\Omega\xi_{3}}{\mid\xi\mid}t)\hat{M}(\xi),
\]
where $\mathcal{I}$ is the identity matrix of order $3\times3$ and $\hat
{M}(\xi)$ is the matrix
\[
\hat{M}(\xi)=\frac{1}{\mid\xi\mid}\left(
\begin{matrix}
0 & \xi_{3} & -\xi_{2}\\
-\xi_{3} & 0 & \xi_{1}\\
\xi_{2} & -\xi_{1} & 0
\end{matrix}
\right)  ,
\]
give us similar estimates in FBM-spaces. Moreover, the fact that functions
$\cos$ and $\sin$ are bounded, and that $\parallel\hat{M}(\xi)\parallel\leq2$
for all $\xi\in\mathbb{R}^{3}$, yield estimates independent of the Coriolis
parameter $\Omega$.

Thus, one also can prove the following result for the fractional
Navier-Stokes-Coriolis system (FNSC), that is, system (\ref{FSBC}) with
$\theta\equiv0$ and $\mathcal{N}=0$. Let $I=(0,\infty)$, $0\leq\mu<3$, $1\leq
q<\infty$ and $s=4-2\alpha-\frac{3-\mu}{q}$. Let $\Omega\in\mathbb{R}$ and
$u_{0}\in\mathcal{FN}_{q,\mu,r}^{s}$ be a divergence-free vector f{}ield.
Assume the conditions (i), (ii) and (iii) in Theorem \ref{TH1}. There are two
constants $\varepsilon=\varepsilon(\nu,\alpha)>0$ and $C=C(\nu,\alpha)>0$
(independent of $\Omega$) such that (FNSC) has a unique global mild solution
$u$ in the closed ball $\{w\in\mathcal{X}_{r};\parallel w\parallel_{X_{r}}\leq
C\varepsilon\}$ satisfying $\nabla\cdot u=0$ provided that $\parallel
u_{0}\parallel_{\mathcal{FN}_{q,\mu,r}^{s}}\leq\varepsilon$. In particular,
considering the critical dissipation $\alpha=1/2$, we have the well-posedness
in critical Fourier-Besov-Morrey spaces for the Navier-Stokes system with
critical dissipation, that is, (FNSC) with the Coriolis parameter $\Omega=0$
and $\alpha=1/2$, which extends the corresponding result of \cite{Xiao-CW2014}
in Fourier-Besov spaces.
\end{obs}

\section{Proofs}

\label{Proof} In this section we prove Theorem \ref{TH1}. For that, we shall
establish some estimates for the fractional Boussinesq-Coriolis-Stratif{}%
ication semigroup (\ref{Semigroup}) and the bilinear term (\ref{bilinearterm}).

\subsection{Linear and Bilinear Estimates}

\begin{lema}
Let $I=(0,\infty)$, $0\leq\mu<3$, $1\leq q_{2}\leq q_{1}<\infty$, $1\leq
r\leq\infty$ and $s\in\mathbb{R}$. For each pair $(\Omega,\mathcal{N}%
)\in(\mathbb{R}-\{0\})^{2},$ consider $L=\max\{2,\frac{\mid\Omega\mid
}{\mathcal{N}\sqrt{g}},\frac{\mathcal{N}\sqrt{g}}{\mid\Omega\mid}\}$. Then,
there is a constant $C>0$, independent of the parameters $\Omega$ and
$\mathcal{N}$, such that
\begin{equation}
\parallel S_{\Omega,N}^{\alpha}(t)v_{0}\parallel_{\mathcal{FN}_{q_{2},\mu
,r}^{s}}\leq CL\nu^{-\frac{1}{2\alpha}(\frac{3-\mu}{q_{2}}-\frac{3-\mu}{q_{1}%
})}t^{-\frac{1}{2\alpha}(\frac{3-\mu}{q_{2}}-\frac{3-\mu}{q_{1}})}\parallel
v_{0}\parallel_{\mathcal{FN}_{q_{1},\mu,r}^{s}}, \label{estsemi1}%
\end{equation}%
\begin{equation}
\parallel S_{\Omega,N}^{\alpha}(t)v_{0}\parallel_{\mathcal{L}^{\infty
}(I;\mathcal{FN}_{q_{1},\mu,r}^{s})}\leq CL\parallel v_{0}\parallel
_{\mathcal{FN}_{q_{1},\mu,r}^{s}}\,\,\mbox{and} \label{estsemi2}%
\end{equation}%
\begin{equation}
\parallel S_{\Omega,N}^{\alpha}(t)v_{0}\parallel_{\mathcal{L}^{1}%
(I;\mathcal{FN}_{q_{1},\mu,r}^{s+2\alpha})}\leq\frac{CL}{\nu}\parallel
v_{0}\parallel_{\mathcal{FN}_{q_{1},\mu,r}^{s}}, \label{estsemi3}%
\end{equation}
for all $v_{0}\in\mathcal{FN}_{q_{1},\mu,r}^{s}$.
\end{lema}

\textbf{Proof.} Consider $p\geq q_{1}^{\prime}\geq1$ such that $\frac{1}%
{q_{2}}=\frac{1}{p}+\frac{1}{q_{1}}$. From the expression of the semigroup
(\ref{Semigroup}), the bound (\ref{EstMatrixM}) of the matrix $M_{j}$ for
$j=1,2,3$, H\"{o}lder's inequality and the scaling property $\parallel
f(\lambda\cdot)\parallel_{p,\mu}=\lambda^{-\frac{n-\mu}{p}}\parallel
f\parallel_{p,\mu}$, it follows that
\begin{align*}
\parallel\varphi_{j}[S_{\Omega,N}^{\alpha}(t)]\symbol{94}\hat{v}_{0}%
\parallel_{q_{2},\mu}  &  =\parallel\varphi_{j}e^{-\nu t\mid\xi\mid^{2\alpha}%
}[\cos(\frac{\mid\xi\mid^{\prime}}{\mid\xi\mid}t)M_{1}(\xi)\hat{v}_{0}%
(\xi)+\sin(\frac{\mid\xi\mid^{\prime}}{\mid\xi\mid}t)M_{2}(\xi)\hat{v}_{0}%
(\xi)+M_{3}(\xi)\hat{v}_{0}(\xi)]\parallel_{q_{2},\mu}\\
&  \leq CL\parallel e^{-\nu t\mid\xi\mid^{2\alpha}}\varphi_{j}\hat{v}_{0}%
(\xi)\parallel_{q_{2},\mu}\\
&  =CL\parallel e^{-\mid(\nu t)^{\frac{1}{2\alpha}}\xi\mid^{2\alpha}}%
\varphi_{j}\hat{v}_{0}(\xi)\parallel_{q_{2},\mu}\\
&  \leq CL\parallel e^{-\mid(\nu t)^{\frac{1}{2\alpha}}\xi\mid^{2\alpha}%
}\parallel_{p,\mu}\parallel\varphi_{j}\hat{v}_{0}\parallel_{q_{1},\mu}\\
&  \leq CL\,(\nu t)^{-\frac{1}{2\alpha}(\frac{3-\mu}{p})}\parallel e^{-\mid
\xi\mid^{2\alpha}}\parallel_{p,\mu}\parallel\varphi_{j}\hat{v}_{0}%
\parallel_{q_{1},\mu}\\
&  \leq CL\,\nu^{-\frac{1}{2\alpha}(\frac{3-\mu}{q_{2}}-\frac{3-\mu}{q_{1}}%
)}t^{-\frac{1}{2\alpha}(\frac{3-\mu}{q_{2}}-\frac{3-\mu}{q_{1}})}%
\parallel\varphi_{j}\hat{v}_{0}\parallel_{q_{1},\mu},
\end{align*}
where we used that $\parallel e^{-\mid\xi\mid^{2\alpha}}\parallel_{p,\mu
}<\infty$. Multiplying by $2^{js}$ and then applying the $l^{r}(\mathbb{Z}%
)$-norm on both sides in the above inequality we get (\ref{estsemi1}). Observe
that, when $q_{1}=q_{2}$, we have
\[
2^{js}\parallel\varphi_{j}[S_{\Omega,N}^{\alpha}(t)]\symbol{94}\hat{u}%
_{0}\parallel_{q_{1},\mu}\leq CL\,2^{js}\parallel\varphi_{j}\hat{u}%
_{0}\parallel_{q_{1},\mu}%
\]
and thus, taking the $L^{\infty}$-norm on $I$ and then the $l^{r}(\mathbb{Z}%
)$-norm, we obtain (\ref{estsemi2}). The last inequality (\ref{estsemi3})
follows from the estimates
\begin{align}
\parallel\parallel\varphi_{j}[S_{\Omega,N}^{\alpha}(t)]\symbol{94}\hat{u}%
_{0}\parallel_{q_{1},\mu}\parallel_{L^{1}(I)}  &  \leq CL\parallel\parallel
e^{-\nu t\mid\xi\mid^{2\alpha}}\varphi_{j}\hat{u}_{0}\parallel_{q_{1},\mu
}\parallel_{L^{1}(I)}\nonumber\\
&  \leq CL\parallel e^{-\nu t2^{2\alpha(j-1)}}\parallel_{L^{1}(I)}%
\parallel\varphi_{j}\hat{u}_{0}\parallel_{q_{1},\mu}\nonumber\\
&  \leq\frac{CL}{\nu}2^{-2\alpha j}\parallel\varphi_{j}\hat{u}_{0}%
\parallel_{q_{1},\mu}, \label{aux-lin-1}%
\end{align}
after multiplying by $2^{j(s+2\alpha)}$ and then applying the $l^{r}%
(\mathbb{Z})$-norm on both sides of (\ref{aux-lin-1}). \fin

Now we shall establish some estimates for the linear operator $\zeta
_{\Omega,N}^{\alpha}$, which is linked to (\ref{bilinearterm}) and def{}ined
in Fourier variables by
\[
\lbrack\zeta_{\Omega,N}^{\alpha}(f)]\symbol{94}(\xi,t)=\displaystyle{\int
_{0}^{t}}[S_{\Omega,N}^{\alpha}(t-\tau)]\symbol{94}(\xi)[\tilde{\mathbb{P}%
}]\symbol{94}(\xi)\hat{f}(\xi,\tau)d\,\tau.
\]

\begin{lema}
\label{LemmaZeta} Let $I=(0,\infty)$, $1\leq r\leq\infty$, $0\leq\mu<3$,
$1\leq q<\infty$ and $s=4-2\alpha-\frac{3-\mu}{q}$. Assume that $(\Omega
,\mathcal{N})\in(\mathbb{R}-\{0\})^{2}$ and consider $L=\max\{2,\frac
{\mid\Omega\mid}{\mathcal{N}\sqrt{g}},\frac{\mathcal{N}\sqrt{g}}{\mid
\Omega\mid}\}.$ There is a constant $C=C(\alpha)>0$ such that
\begin{equation}
\parallel\zeta_{\Omega,N}^{\alpha}(f)\parallel_{\mathcal{L}^{\infty
}(I;\mathcal{FN}_{q,\mu,r}^{s})}\leq CL\parallel f\parallel_{\mathcal{L}%
^{1}(I;\mathcal{FN}_{q,\mu,r}^{s})}, \label{LL1}%
\end{equation}%
\begin{equation}
\parallel\zeta_{\Omega,N}^{\alpha}(f)\parallel_{\mathcal{L}^{\infty
}(I;\mathcal{FN}_{q,\mu,r}^{s+2\alpha})}\leq\frac{CL}{\nu}\parallel
f\parallel_{\mathcal{L}^{\infty}(I;\mathcal{FN}_{q,\mu,r}^{s})}, \label{LL2}%
\end{equation}%
\begin{equation}
\parallel\zeta_{\Omega,N}^{\alpha}(f)\parallel_{\mathcal{L}^{1}(I;\mathcal{FN}%
_{q,\mu,r}^{s+2\alpha})}\leq\frac{CL}{\nu}\parallel f\parallel_{\mathcal{L}%
^{1}(I;\mathcal{FN}_{q,\mu,r}^{s})}, \label{LL3}%
\end{equation}
for all $f\in\mathcal{L}^{\infty}(I;\mathcal{FN}_{q,\mu,r}^{s})$ or
$f\in\mathcal{L}^{1}(I;\mathcal{FN}_{q,\mu,r}^{s}),$ according to the
corresponding estimate.
\end{lema}

\textbf{Proof.} Using $supp\,(\varphi_{j})\subset D_{j}=\{\xi\in\mathbb{R}%
^{3};\,2^{j-1}\leq\mid\xi\mid\leq2^{j+2}\}$, $\parallel\lbrack\tilde
{\mathbb{P}}]\symbol{94}(\xi)\parallel\leq2$, and the fact that $L$ is a bound
of the matrix $M_{j}$ for $j=1,2,3$, we can use Young's inequality in the
time-variable in order to estimate
\begin{align*}
\parallel\parallel\varphi_{j}[\zeta_{\Omega,N}^{\alpha}(f)]\symbol{94}%
(\cdot,t)\parallel_{q,\mu}\parallel_{L^{\infty}(I)}  &  \leq\parallel
\displaystyle{\int_{0}^{t}}\parallel\mid\lbrack S_{\Omega,N}^{\alpha}%
(t-\tau)]\symbol{94}(\xi)\mid\varphi_{j}\hat{f}\parallel_{q,\mu}%
\parallel_{L^{\infty}(I)}\\
&  \leq CL\parallel\displaystyle{\int_{0}^{t}}\parallel e^{-\nu(t-\tau)\mid
\xi\mid^{2\alpha}}\varphi_{j}\hat{f}\parallel_{q,\mu}\parallel_{L^{\infty}%
(I)}\\
&  \leq CL\parallel\displaystyle{\int_{0}^{t}}e^{-\nu(t-\tau)2^{2\alpha(j-1)}%
}\parallel\varphi_{j}\hat{f}\parallel_{q,\mu}\parallel_{L^{\infty}(I)}\\
&  \leq CL\parallel e^{-\nu t2^{2\alpha(j-1)}}\parallel_{L^{\infty}%
(I)}\parallel\varphi_{j}\hat{f}\parallel_{L^{1}(I;\mathcal{M}_{q,\mu})}\\
&  \leq CL\parallel\varphi_{j}\hat{f}\parallel_{L^{1}(I;\mathcal{M}_{q,\mu})}.
\end{align*}
Multiplying by $2^{js}$ and then applying the $l^{r}(\mathbb{Z})$-norm on both
sides of the above inequality, we get (\ref{LL1}). For (\ref{LL2}), we use
Young's inequality in the time-variable in order to obtain
\begin{align}
\parallel\parallel\varphi_{j}[\zeta_{\Omega,N}^{\alpha}]\symbol{94}%
(\cdot,t)\parallel_{q,\mu}\parallel_{L^{\infty}(I)}  &  \leq CL\parallel
\displaystyle{\int_{o}^{t}}e^{-\nu(t-\tau)2^{2\alpha(j-1)}}\parallel
\varphi_{j}\hat{f}(\cdot,\tau)\parallel_{q,\mu}d\,\tau\parallel_{L^{\infty
}(I)}\nonumber\\
&  \leq CL\parallel e^{-\nu t2^{2\alpha(j-1)}}\parallel_{L^{1}(I)}%
\parallel\varphi_{j}\hat{f}\parallel_{L^{\infty}(I;\mathcal{M}_{q,\mu}%
)}\nonumber\\
&  \leq\frac{CL}{\nu}\,2^{-2\alpha j}\parallel\varphi_{j}\hat{f}%
\parallel_{L^{\infty}(I;\mathcal{M}_{q,\mu})}, \label{aux-Lin-2}%
\end{align}
which yields (\ref{LL2}) after multiplying by $2^{j(s+2\alpha)}$ and then
applying the $l^{r}(\mathbb{Z})$-norm on its both sides. For the last
estimate, we have that
\begin{align}
\parallel\parallel\varphi_{j}[\zeta_{\Omega,N}^{\alpha}]\symbol{94}%
(\cdot,t)\parallel_{q,\mu}\parallel_{L^{1}(I)}  &  \leq CL\parallel
\displaystyle{\int_{o}^{t}}e^{-\nu(t-\tau)2^{2\alpha(j-1)}}\parallel
\varphi_{j}\hat{f}(\cdot,\tau)\parallel_{q,\mu}d\,\tau\parallel_{L^{1}%
(I)}\nonumber\\
&  \leq CL\parallel e^{-\nu t2^{2\alpha(j-1)}}\parallel_{L^{1}(I)}%
\parallel\varphi_{j}\hat{f}\parallel_{L^{1}(I;\mathcal{M}_{q,\mu})}\nonumber\\
&  \leq\frac{CL}{\nu}\,2^{-2\alpha j}\parallel\varphi_{j}\hat{f}%
\parallel_{L^{1}(I;\mathcal{M}_{q,\mu})}. \label{aux-Lin-33}%
\end{align}
Now it is sufficient to multiply by $2^{(s+2\alpha)j}$ and then take the
$l^{r}(\mathbb{Z})$-norm on both sides of (\ref{aux-Lin-33}) to obtain
(\ref{LL3}). \fin

\begin{lema}
\label{LemmaBilinear} Let $I=(0,\infty)$, $0\leq\mu<3$, $1\leq q<\infty$,
$1\leq r\leq\infty$ and $s=4-2\alpha-\frac{3-\mu}{q}$. Assume that $\alpha$
and $r$ satisfy either the following conditions

\begin{itemize}
\item[(i)] $\frac{1}{2}<\alpha<\frac{5}{2}-\frac{3-\mu}{2q}$ and $1\leq
r\leq\infty$, or

\item[(ii)] $\frac{1}{2}<\alpha=\frac{5}{2}-\frac{3-\mu}{2q}$, $\mu=0$ and
$1\leq q\leq r\leq2$, or

\item[(iii)] $\alpha=\frac{1}{2}$ and $r=1$.
\end{itemize}

Let $(\Omega,\mathcal{N})\in(\mathbb{R}-\{0\})^{2}$ and $L=\max\{2,\frac
{\mid\Omega\mid}{\mathcal{N}\sqrt{g}},\frac{\mathcal{N}\sqrt{g}}{\mid
\Omega\mid}\}.$ Then, there is a constant $C=C(\alpha)>0$ such that
\begin{equation}
\parallel B(v,w)\parallel_{\mathcal{X}_{r}}\leq CL\max\{1,\frac{1}{\nu
}\}\parallel v\parallel_{\mathcal{X}_{r}}\parallel w\parallel_{\mathcal{X}%
_{r}}, \label{Bilinearbounded}%
\end{equation}
for all $v,w\in\mathcal{X}_{r}=\mathcal{L}^{\infty}(I;\mathcal{FN}_{q,\mu
,r}^{s})\cap\mathcal{L}^{1}(I;\mathcal{FN}_{q,\mu,r}^{s+2\alpha})$.
\end{lema}

\textbf{Proof. }First we shall establish estimate (\ref{Bilinearbounded}) for
the assumption $(i)$. Since $supp\,\varphi_{j}\subset\{\xi\in\mathbb{R}%
^{3};\mid\xi\mid\leq\frac{4}{3}2^{j+1}\}$ and$\parallel\lbrack\tilde
{\mathbb{P}}]\symbol{94}(\xi)\parallel\leq2,$ Bernstein-type inequality
yields
\[
\parallel\varphi_{j}[i\xi,0]^{T}[v\otimes w]\symbol{94}(\xi,t)\parallel
_{q,\mu}\leq C2^{j+1}\parallel\varphi_{j}[v\otimes w]\symbol{94}%
(\xi,t)\parallel_{q,\mu}.
\]
Now, we use Lemma \ref{LemmaZeta} with $f=\tilde{\nabla}\cdot(v\otimes w)$ to
obtain
\begin{align*}
\parallel B(v,w)\parallel_{\mathcal{L}^{1}(I;\mathcal{FN}_{q,\mu,r}%
^{s+2\alpha})}+\parallel B(v,w)\parallel_{\mathcal{L}^{\infty}(I;\mathcal{FN}%
_{q,\mu,r}^{s})}  &  \leq CL\,\max\{1,\frac{1}{\nu}\}\parallel f\parallel
_{\mathcal{L}^{1}(I;\mathcal{FN}_{q,\mu,r}^{s})}\\
&  \leq CL\,\max\{1,\frac{1}{\nu}\}\parallel v\otimes w\parallel
_{\mathcal{L}^{1}(I;\mathcal{FN}_{q,\mu,r}^{s+1})}.
\end{align*}
To conclude we are going to show that
\[
\parallel v\otimes w\parallel_{\mathcal{L}^{1}(I;\mathcal{FN}_{q,\mu,r}%
^{s+1})}\leq C\parallel v\parallel_{\mathcal{X}_{r}}\parallel w\parallel
_{\mathcal{X}_{r}}.
\]
For that, we consider Bony's decomposition
\begin{align*}
\varphi_{j}[vw]\symbol{94}  &  =\displaystyle{\sum_{\mid k-j\mid\leq4}}%
\varphi_{j}[(S_{k-1}v)\symbol{94}\ast(\varphi_{k}\hat{w})]+\displaystyle{\sum
_{\mid k-j\mid\leq4}}\varphi_{j}[(S_{k-1}w)\symbol{94}\ast(\varphi_{k}\hat
{v})]+\displaystyle{\sum_{k\geq j-2}}\varphi_{j}[(\varphi_{k}\hat{v}%
)\ast(\tilde{\varphi}_{k}\hat{w})]\\
&  =I_{1}+I_{2}+I_{3}.
\end{align*}
Estimate for the term $I_{1}$: First note that $\parallel\varphi_{k}\hat
{u}\parallel_{L^{1}}\leq C2^{k(3-\frac{3-\mu}{q})}\parallel\varphi_{k}\hat
{u}\parallel_{q,\mu}$follows from the Bernstein-type inequality with $\beta
=0$, $(p_{2},\mu_{2})=(1,0)$ and $(p_{1},\mu_{1})=(q,\mu)$. Next, we proceed
as follows:
\begin{align}
\parallel I_{1}\parallel_{L^{1}(I;\mathcal{M}_{q,\mu})}  &  \leq
\,\displaystyle{\int_{I}}\displaystyle{\sum_{\mid k-j\mid\leq4}}%
\parallel\varphi_{j}[(S_{k-1}v)\symbol{94}\ast(\varphi_{k}\hat{w}%
)]\parallel_{q,\mu}d\,t\nonumber\\
&  \leq C\displaystyle{\int_{I}}\displaystyle{\sum_{\mid k-j\mid\leq4}%
}\displaystyle{(\displaystyle{\sum_{k^{\prime}<k-2}}\parallel\varphi
_{k^{\prime}}\hat{v}\parallel_{L^{1}})}\parallel\varphi_{k}\hat{w}%
\parallel_{q,\mu}d\,t\nonumber\\
&  \leq C\displaystyle{\int_{I}}\displaystyle{\sum_{\mid k-j\mid\leq4}%
}\,\displaystyle{\sum_{k^{\prime}<k-2}}2^{k^{\prime}(-1+2\alpha)}2^{k^{\prime
}(4-2\alpha-\frac{3-\mu}{q})}\parallel\varphi_{k^{\prime}}\hat{v}%
\parallel_{q,\mu}\parallel\varphi_{k}\hat{w}\parallel_{q,\mu}d\,t\nonumber\\
&  \leq C\displaystyle{\int_{I}}\displaystyle{\sum_{\mid k-j\mid\leq4}%
}2^{k(-1+2\alpha)}\,\parallel2^{k^{\prime}(4-2\alpha-\frac{3-\mu}{q}%
)}\parallel\varphi_{k^{\prime}}\hat{v}\parallel_{L^{\infty}(I;\mathcal{M}%
_{q,\mu})}\parallel_{l^{r}(\mathbb{Z})}\parallel\varphi_{k}\hat{w}%
\parallel_{q,\mu}d\,t\nonumber\\
&  \leq C\parallel v\parallel_{\mathcal{L}^{\infty}(I;\mathcal{FN}_{q,\mu
,r}^{s})}\displaystyle{\int_{I}}\displaystyle{\sum_{\mid k-j\mid\leq4}%
}2^{k(-1+2\alpha)}\parallel\varphi_{k}\hat{w}\parallel_{q,\mu}d\,t\nonumber\\
&  \leq C\parallel v\parallel_{\mathcal{L}^{\infty}(I;\mathcal{FN}_{q,\mu
,r}^{s})}2^{j(-5+2\alpha+\frac{3-\mu}{q})}\,\displaystyle{\sum_{k\in
\mathbb{Z}}}2^{-(j-k)(-5+2\alpha+\frac{3-\mu}{q})}\chi_{\{l;\mid l\mid
\leq4\}(j-k)}\times\nonumber\\
&  2^{k(4-\frac{3-\mu}{q})}\parallel\varphi_{k}\hat{w}\parallel_{L^{1}%
(I;\mathcal{M}_{q,\mu})}\nonumber\\
&  \leq C\parallel v\parallel_{\mathcal{L}^{\infty}(I;\mathcal{FN}_{q,\mu
,r}^{s})}2^{j(-5+2\alpha+\frac{3-\mu}{q})}(a_{l}\ast b_{k})_{j},
\label{aux-Lin4}%
\end{align}
where $a_{l}=2^{-l(-5+2\alpha+\frac{3-\mu}{q})}\chi_{\{l;\mid l\mid\leq4\}}$
and $b_{k}=2^{k(4-\frac{3-\mu}{q})}\parallel\varphi_{k}\hat{w}\parallel
_{L^{1}(I;\mathcal{M}_{q,\mu})}$. Multiplying by $2^{j(5-2\alpha-\frac{3-\mu
}{q})}$, taking the $l^{r}(\mathbb{Z})$-norm on both sides of (\ref{aux-Lin4}%
), and using Young's inequality for series, we get
\begin{align*}
\parallel2^{j(5-2\alpha-\frac{3-\mu}{q})}\parallel I_{1}\parallel
_{L^{1}(I;\mathcal{M}_{q,\mu})}\parallel_{l^{r}(\mathbb{Z})}  &  \leq
C\parallel v\parallel_{\mathcal{L}^{\infty}(I;\mathcal{FN}_{q,\mu,r}^{s}%
)}\parallel a_{l}\parallel_{l^{1}}\parallel b_{k}\parallel_{l^{r}}\\
&  \leq C\parallel v\parallel_{\mathcal{L}^{\infty}(I;\mathcal{FN}_{q,\mu
,r}^{s})}\parallel w\parallel_{\mathcal{L}^{1}(I;\mathcal{FN}_{q,\mu
,r}^{s+2\alpha})}%
\end{align*}
and then
\begin{equation}
\parallel2^{j(s+1)}\parallel I_{1}\parallel_{L^{1}(I;\mathcal{M}_{q,\mu}%
)}\parallel_{l^{r}(\mathbb{Z})}\leq C\parallel v\parallel_{\mathcal{X}_{r}%
}\parallel w\parallel_{\mathcal{X}_{r}}. \label{A1}%
\end{equation}
Similar computations lead us
\begin{equation}
\parallel2^{j(s+1)}\parallel I_{2}\parallel_{L^{1}(I;\mathcal{M}_{q,\mu}%
)}\parallel_{l^{r}(\mathbb{Z})}\leq C\parallel v\parallel_{\mathcal{X}_{r}%
}\parallel w\parallel_{\mathcal{X}_{r}}. \label{A2}%
\end{equation}
Finally we estimate the parcel $I_{3}$:
\begin{align*}
\parallel I_{3}\parallel_{L^{1}(I;\mathcal{M}_{q,\mu})}  &  \leq
\,\displaystyle{\int_{I}}\displaystyle{\sum_{k\geq j-2}}\parallel\varphi
_{j}[(\varphi_{k}\hat{v})\ast(\tilde{\varphi}_{k}\hat{w})]\parallel_{q,\mu
}d\,t\\
&  \leq C\,\displaystyle{\int_{I}}\displaystyle{\sum_{k\geq j-2}}%
\parallel\varphi_{k}\hat{v}\parallel_{q,\mu}\parallel\tilde{\varphi}_{k}%
\hat{w}\parallel_{L^{1}}d\,t\\
&  \leq C\,\displaystyle{\int_{I}}\displaystyle{\sum_{k\geq j-2}}%
\parallel\varphi_{k}\hat{v}\parallel_{q,\mu}\displaystyle{\sum_{\mid
k-k^{\prime}\mid\leq1}}\parallel\varphi_{k^{\prime}}\hat{w}\parallel_{L^{1}%
}d\,t\\
&  \leq C\,\displaystyle{\sum_{k\geq j-2}}\displaystyle{\int_{I}%
}\displaystyle{\sum_{\mid k-k^{\prime}\mid\leq1}}2^{k^{\prime}(3-\frac{3-\mu
}{q})}\parallel\varphi_{k^{\prime}}\hat{w}\parallel_{q,\mu}\parallel
\varphi_{k}\hat{v}\parallel_{q,\mu}d\,t\\
&  \leq C\,\displaystyle{\sum_{k\geq j-2}}\displaystyle{\int_{I}%
}\displaystyle{(\displaystyle{\sum_{\mid k-k^{\prime}\mid\leq1}}2^{k^{\prime
}(-1+2\alpha)}2^{k^{\prime}(4-2\alpha-\frac{3-\mu}{q})}\parallel
\varphi_{k^{\prime}}\hat{w}\parallel_{L^{\infty}(I;\mathcal{M}_{q,\mu})}%
)}\parallel\varphi_{k}\hat{v}\parallel_{q,\mu}d\,t\\
&  \leq C\,\displaystyle{\sum_{k\geq j-2}}\displaystyle{\int_{I}%
}2^{k(-1+2\alpha)}\parallel2^{k^{\prime}(4-2\alpha-\frac{3-\mu}{q}%
)}\,\parallel\varphi_{k^{\prime}}\hat{w}\parallel_{L^{\infty}(I;\mathcal{M}%
_{q,\mu})}\parallel_{l^{r}(\mathbb{Z})}\parallel\varphi_{k}\hat{v}%
\parallel_{q,\mu}d\,t\\
&  \leq C\,\parallel w\parallel_{\mathcal{L}^{\infty}(I;\mathcal{FN}_{q,\mu
,r}^{s})}\displaystyle{\sum_{k\geq j-2}}2^{k(-1+2\alpha)}\parallel\varphi
_{k}\hat{v}\parallel_{L^{1}(I;\mathcal{M}_{q,\mu})}\\
&  \leq C\,\parallel w\parallel_{\mathcal{L}^{\infty}(I;\mathcal{FN}_{q,\mu
,r}^{s})}\displaystyle{\sum_{k\geq j-2}}2^{k(-5+2\alpha+\frac{3-\mu}{q}%
)}2^{k(4-\frac{3-\mu}{q})}\parallel\varphi_{k}\hat{v}\parallel_{L^{1}%
(I;\mathcal{M}_{q,\mu})}\\
&  =C\,\parallel w\parallel_{\mathcal{L}^{\infty}(I;\mathcal{FN}_{q,\mu,r}%
^{s})}2^{j(-5+2\alpha+\frac{3-\mu}{q})}\displaystyle{\sum_{k\in\mathbb{Z}}%
}2^{-(j-k)(-5+2\alpha+\frac{3-\mu}{q})}\chi_{\{l;l\leq2\}}(j-k)\times\\
&  2^{k(4-\frac{3-\mu}{q})}\parallel\varphi_{k}\hat{v}\parallel_{L^{1}%
(I;\mathcal{M}_{q,\mu})}.
\end{align*}
Thus
\[
2^{j(s+1)}\parallel I_{3}\parallel_{L^{1}(I;\mathcal{M}_{q,\mu})}\leq
C\,\parallel w\parallel_{\mathcal{L}^{\infty}(I;\mathcal{FN}_{q,\mu,r}^{s}%
)}(a_{l}\ast b_{k})_{j},
\]
with $a_{l}=2^{-l(-5+2\alpha+\frac{3-\mu}{q})}\chi_{\{l;l\leq2\}}$ and
$b_{k}=2^{k(4-\frac{3-\mu}{q})}\parallel\varphi_{k}\hat{v}\parallel
_{L^{1}(I;\mathcal{M}_{q,\mu})}$. After applying the $l^{r}(\mathbb{Z})$-norm
on both sides of the above inequality, and using Young's inequality for
series, we arrive at
\begin{align*}
\parallel2^{j(s+1)}\parallel I_{3}\parallel_{L^{1}(I;\mathcal{M}_{q,\mu}%
)}\parallel_{l^{r}(\mathbb{Z})}  &  \leq C\parallel w\parallel_{\mathcal{L}%
^{\infty}(I;\mathcal{FN}_{q,\mu,r}^{s})}\parallel a_{l}\parallel
_{l^{1}(\mathbb{Z})}\parallel b_{k}\parallel_{l^{r}(\mathbb{Z})}\\
&  \leq C\parallel v\parallel_{\mathcal{L}^{1}(I;\mathcal{FN}_{q,\mu
,r}^{s+2\alpha})}\parallel w\parallel_{\mathcal{L}^{\infty}(I;\mathcal{FN}%
_{q,\mu,r}^{s})}\\
&  \leq\parallel v\parallel_{\mathcal{X}_{r}}\parallel w\parallel
_{\mathcal{X}_{r}},
\end{align*}
where we used the assumption $-5+2\alpha+\frac{3-\mu}{q}<0$, i.e.,
$\alpha<\frac{5}{2}-\frac{3-\mu}{2q}$. Therefore, using (\ref{A1}), (\ref{A2})
and the last estimate, we obtain
\begin{align*}
\parallel v\otimes w\parallel_{\mathcal{L}^{1}(I;\mathcal{FN}_{q,\mu,r}%
^{s+1})}  &  \leq\parallel2^{j(s+1)}\parallel I_{1}\parallel_{L^{1}%
(I;\mathcal{M}_{q,\mu})}\parallel_{l^{r}(\mathbb{Z})}+\parallel2^{j(s+1)}%
\parallel I_{2}\parallel_{L^{1}(I;\mathcal{M}_{q,\mu})}\parallel
_{l^{r}(\mathbb{Z})}+\\
&  \parallel2^{j(s+1)}\parallel I_{3}\parallel_{L^{1}(I;\mathcal{M}_{q,\mu}%
)}\parallel_{l^{r}(\mathbb{Z})}\\
&  \leq C\parallel v\parallel_{\mathcal{X}_{r}}\parallel w\parallel
_{\mathcal{X}_{r}}.
\end{align*}
These computations give us (\ref{Bilinearbounded}) under the conditions of
item (i).\newline

For the conditions in item (ii), we can proceed similarly to the proof of
\cite[ Theorem 5]{Xiao-CW2014} in order to obtain (\ref{Bilinearbounded}). The
details are left to the reader.\newline

Finally, assume the conditions in item (iii). For the f{}irst term $I_{1}$, we
have
\begin{align*}
&  \parallel I_{1}\parallel_{L^{1}(I;\mathcal{M}_{q,\mu})}\leq
C\,\displaystyle{\sum_{\mid k-j\mid\leq4}}\,\displaystyle{\sum_{k^{\prime
}<k-2}}2^{k^{\prime}(3-\frac{3-\mu}{q})}\parallel\varphi_{k^{\prime}}\hat
{v}\parallel_{L^{\infty}(I;\mathcal{M}_{q,\mu})}\parallel\varphi_{k}\hat
{w}\parallel_{L^{1}(I;\mathcal{M}_{q,\mu})}\\
&  \leq C\,\parallel v\parallel_{\mathcal{L}^{\infty}(I;\mathcal{FN}_{q,\mu
,1}^{3-\frac{3-\mu}{q}})}2^{-j(4-\frac{3-\mu}{q})}(a_{l}\ast b_{k})_{j}\text{
},
\end{align*}
where $a_{l}=2^{l(4-\frac{3-\mu}{q})}\chi_{\{l;\mid l\mid\leq4\}}$ and
$b_{k}=2^{k(4-\frac{3-\mu}{q})}\parallel\varphi_{k}\hat{w}\parallel
_{L^{1}(I;\mathcal{M}_{q,\mu})}$. Multiplying by $2^{j(4-\frac{3-\mu}{q})}$
and then taking the $l^{1}$-norm on both sides of the above estimate, and
using Young's inequality for series, we get
\[
\parallel2^{j(4-\frac{3-\mu}{q})}\parallel I_{1}\parallel_{L^{1}%
(I;\mathcal{M}_{q,\mu})}\parallel_{l^{1}}\leq C\,\parallel v\parallel
_{\mathcal{L}^{\infty}(I;\mathcal{FN}_{q,\mu,1}^{3-\frac{3-\mu}{q}})}\parallel
w\parallel_{\mathcal{L}^{1}(I;\mathcal{FN}_{q,\mu,1}^{4-\frac{3-\mu}{q}})}.
\]
Similar computations lead us to
\[
\parallel2^{j(4-\frac{3-\mu}{q})}\parallel I_{2}\parallel_{L^{1}%
(I;\mathcal{M}_{q,\mu})}\parallel_{l^{1}}\leq C\,\parallel w\parallel
_{\mathcal{L}^{\infty}(I;\mathcal{FN}_{q,\mu,1}^{3-\frac{3-\mu}{q}})}\parallel
v\parallel_{\mathcal{L}^{1}(I;\mathcal{FN}_{q,\mu,1}^{4-\frac{3-\mu}{q}})}.
\]
\newline Now we consider the estimate for the term $I_{3}$. We have that
\begin{align*}
&  \parallel I_{3}\parallel_{L^{1}(I;\mathcal{M}_{q,\mu})}\leq
C\,\displaystyle{\sum_{k\geq j-2}}\,\displaystyle{\int_{I}}%
\,\displaystyle{\sum_{\mid k-k^{\prime}\mid\leq1}}2^{k^{\prime}(3-\frac{3-\mu
}{q})}\parallel\varphi_{k^{\prime}}\hat{w}\parallel_{q,\mu}\parallel
\varphi_{k}\hat{v}\parallel_{q,\mu}\,dt\\
&  \leq C\,\displaystyle{\sum_{k\geq j-2}}\,\displaystyle{\sum_{\mid
k-k^{\prime}\mid\leq1}}2^{k^{\prime}(3-\frac{3-\mu}{q})}\parallel
\varphi_{k^{\prime}}\hat{w}\parallel_{L^{\infty}(I;\mathcal{M}_{q,\mu}%
)}\parallel\varphi_{k}\hat{v}\parallel_{L^{1}(I;\mathcal{M}_{q,\mu})}\\
&  \leq C\,\parallel w\parallel_{\mathcal{L}^{\infty}(I;\mathcal{FN}_{q,\mu
,1}^{3-\frac{3-\mu}{q}})}2^{-j(4-\frac{3-\mu}{q})}(a_{l}\ast b_{k})_{j}\text{
},
\end{align*}
where $a_{l}=2^{l(4-\frac{3-\mu}{q})}\chi_{\{l;l\leq2\}}$ and $b_{k}%
=2^{k(4-\frac{3-\mu}{q})}\parallel\varphi_{k}\hat{v}\parallel_{L^{1}%
(I;\mathcal{M}_{q,\mu})}$. Multiplying by $2^{j(4-\frac{3-\mu}{q})}$, taking
the $l^{1}$-norm on both sides of the above estimate, and using Young's
inequality for series, we obtain
\[
\parallel2^{j(4-\frac{3-\mu}{q})}\parallel I_{3}\parallel_{L^{1}%
(I;\mathcal{M}_{q,\mu})}\parallel_{l^{1}}\leq C\,\parallel w\parallel
_{\mathcal{L}^{\infty}(I;\mathcal{FN}_{q,\mu,1}^{3-\frac{3-\mu}{q}})}\parallel
v\parallel_{\mathcal{L}^{1}(I;\mathcal{FN}_{q,\mu,1}^{4-\frac{3-\mu}{q}})},
\]
and we are done.\fin

\subsection{Proof of Theorem \ref{TH1}}

In this subsection we shall apply Lemma \ref{fixedpoint}. Consider
$y=S_{\Omega,N}^{\alpha}(t)v_{0}$ and $\mathcal{X}_{r}=\mathcal{L}^{\infty
}(I;\mathcal{FN}_{q,\mu,r}^{s})\cap\mathcal{L}^{1}(I;\mathcal{FN}_{q,\mu
,r}^{s+2\alpha})$ with $I=(0,\infty)$. By estimates (\ref{estsemi2}) and
(\ref{estsemi3}) with $q_{1}=q,$ we obtain
\[
\parallel y\parallel_{\mathcal{X}_{r}}\leq CL\max\{1,\frac{1}{\nu}\}\parallel
v_{0}\parallel_{\mathcal{FN}_{q,\mu,r}^{s}},
\]
where $C$ is a constant. Bilinear estimate (\ref{Bilinearbounded}) yields
\[
\parallel B(v,w)\parallel_{\mathcal{X}_{r}}\leq CL\,\max\{1,\frac{1}{\nu
}\}\parallel v\parallel_{\mathcal{X}_{r}}\parallel w\parallel_{\mathcal{X}%
_{r}},
\]
and we can choose the initial data $v_{0}$ satisfying
\[
\parallel v_{0}\parallel_{\mathcal{FN}_{q,\mu,r}^{s}}\leq\varepsilon<\frac
{1}{4(CL\,\max\{1,\frac{1}{\nu}\})^{2}}.
\]
Therefore, Lemma \ref{fixedpoint} implies that system (\ref{FSBC}) has a
unique global mild solution $v\in\mathcal{X}_{r}$ such that $\parallel
v\parallel_{\mathcal{X}_{r}}\leq2CL\,\max\{1,\frac{1}{\nu}\}\varepsilon$.
Observe that the constant $R$ depends on the parameters $\Omega$ and
$\mathcal{N}$, except when $\mathcal{N}\sqrt{g}/2\leq\Omega\leq2\mathcal{N}%
\sqrt{g}$, in which case the value is equal to $2$. Since $\tilde{\nabla}\cdot
v_{0}=0$ in $\mathcal{S}^{\prime}(\mathcal{R}^{3})$, the mild formulation
(\ref{mildform}) implies that $\tilde{\nabla}\cdot v=0$ for $v=(u^{1}%
,u^{2},u^{3},\sqrt{g}\theta/{\mathcal{N}})$. In order to prove Theorem
\ref{TH1} in the case of items (ii) and (iii), we use the estimate
(\ref{Bilinearbounded}) with the conditions (ii) and (iii) of Lemma
\ref{LemmaBilinear}, respectively. \newline

Finally, we employ the estimate (\ref{fixedpointdependence}) in Lemma
\ref{fixedpoint} in order to show the continuous dependence of $v$ with
respect to the initial data $v_{0}$. In fact, we know that
\begin{align*}
\parallel v-\tilde{v}\parallel_{\mathcal{X}_{r}}  &  \leq\frac{1}%
{1-4(CL\,\max\{1,\frac{1}{\nu}\})^{2}\varepsilon}\parallel S_{\Omega
,N}^{\alpha}(t)v_{0}-S_{\Omega,N}^{\alpha}(t)\tilde{v}_{0}\parallel
_{\mathcal{X}_{r}}\\
&  \leq\frac{CL\,\max\{1,\frac{1}{\nu}\}}{1-4(CL\,\max\{1,\frac{1}{\nu}%
\})^{2}\varepsilon}\parallel v_{0}-\tilde{v}_{0}\parallel_{\mathcal{FN}%
_{q,\mu,r}^{s}}.
\end{align*}
This inequality gives us the Lipschitz continuity for the data-solution map
$v_{0}\mapsto v$ from $\{v_{0}\in\mathcal{FN}_{q,\mu,r}^{s};\parallel
v_{0}\parallel_{\mathcal{FN}_{q,\mu,r}^{s}}\leq\varepsilon\}$ to
$\{v\in\mathcal{X}_{r};\parallel v\parallel_{\mathcal{X}_{r}}\leq
2CL\,\max\{1,\frac{1}{\nu}\}\varepsilon\}$. Finally, the weak time-continuity
statement follows from standard arguments (see, e.g., \cite{Kozo-Yamazaki1994}).\fin

\end{document}